\documentclass[
    ,final            
  ]
  {aipproc}

\layoutstyle{8x11single}

  \usepackage{latexsym}
  \usepackage{amsmath}   
  \usepackage{amsfonts}
  \usepackage{amssymb}
  \usepackage{graphicx}
  \usepackage{amsthm}

  \renewcommand{\mod}{\mathrm{mod}\,\,}
  \newcommand{\alert}{}

  \newcommand{\field}[1]{\mathbb{#1}}
  \newcommand{\C}{\field{C}}
  \newcommand{\R}{\field{R}}

  \newcommand{\BK}{\field{K}}

  \newcommand{\vect}[1]{\ensuremath{\mbox{\textbf{\textit{#1}}}}}
  \newcommand{\svect}[1]{\ensuremath{\mbox{\textbf{\textit{\small #1}}}}}
  \newcommand{\bomega}{\boldsymbol{\omega}}

  \newcommand{\be}{\begin{equation}}
  \newcommand{\ee}{\end{equation}} 
  
  \newcommand{\blue}[1]{{#1}}

  \newcommand{\black}[1]{{#1}}
  \newcommand{\bfr}{\begin{frame}}
  \newcommand{\efr}{\end{frame}}

  \newcommand{\btheta}{\boldsymbol{\theta}}
  \newcommand{\sbtheta}{\boldsymbol{\theta}}
  \newcommand{\sbthetap}{\boldsymbol{\theta^{\prime}}}
  \newcommand{\sbthetaz}{\boldsymbol{\theta_0}}
  
  \newtheorem{theorem}{Theorem}
  \newtheorem{definition}{Definition}
  \newtheorem{lemma}{Lemma}
  \newtheorem{corollary}{Corollary}

\begin{document}

\title{Real Clifford Algebra $Cl_{n,0}, n=2,3(\mod 4)$ Wavelet Transform }

\classification{AMS Subj. Class. 15A66, 42C40, 94A12}
\keywords      {Clifford geometric algebra, Clifford wavelet transform, multidimensional wavelets, continuous wavelets, similitude group}

\author{Eckhard Hitzer}{
  address={Department of Applied Physics, University of Fukui, 910-8507 Japan}
}

\begin{abstract}
We show  how for $n=2,3 (\mod 4)$ continuous Clifford (geometric) algebra (GA) $Cl_n$-valued admissible wavelets can be constructed using the similitude group $SIM(n)$. We strictly aim for real geometric interpretation, and replace the imaginary unit $i \in \C$ therefore with a GA blade squaring to $-1$. Consequences due to non-commutativity arise. We express the admissibility condition in terms of a $Cl_{n}$ Clifford Fourier Transform and then derive a set of important properties such as dilation, translation and rotation covariance, a reproducing kernel, and show how to invert the Clifford wavelet transform. As an example, we introduce Clifford Gabor wavelets. We further invent a generalized Clifford wavelet uncertainty principle.  
\end{abstract}

\maketitle


\section{Multivector functions}

Multivectors
  $M \in Cl_{p,q}, p+q=n,$ have $k$-vector parts ($0\leq k \leq n$):
  \alert{scalar} part
  $Sc(M) = \langle M \rangle = \langle M \rangle_0 = M_0 \in \R$, 
  \alert{vector} part
  $\langle M \rangle_1 \in \R^{p,q}$, 
  \alert{bi-vector} part
  $\langle M \rangle_2$,  \ldots, 
  and
  \alert{pseudoscalar} part $\langle M \rangle_n\in\bigwedge^n\R^{p,q}$
  \begin{equation}\label{eq:MVgrades}
    M  =  \sum_{A} M_{A} \vect{e}_{A}
       =  \langle M \rangle + \langle M \rangle_1 + \langle M \rangle_2 + \ldots +\langle M \rangle_n,
  \end{equation}
  with blade index $A \in \{0,1,2,3,12,23,31,123, \,\ldots \,, 12\ldots n \}$, $M_A\in\R$.
  The \blue{Reverse} of $M \in Cl_{p,q}$ is defined as
  \begin{equation}\label{eq:MVrev}
  \widetilde{M}=\; \sum_{k=0}^{n}(-1)^{\frac{k(k-1)}{2}}\langle M \rangle_k,
  \end{equation}
  it replaces \blue{complex conjugation and quaternion conjugation}.
The \alert{scalar product} of two multivectors $M, \widetilde{N} \in Cl_{p,q}$ is defined as
  \be
    M \ast \widetilde{N} 
    = \langle M\widetilde{N} \rangle 
    = \langle M\widetilde{N} \rangle_0.
  \ee
  For $M, \widetilde{N} \in Cl_{n}=Cl_{n,0}$ we get $M\ast \widetilde{N}=\sum_{A} M_A N_A.$
  The \blue{modulus} $|M|$ of a multivector $M \in Cl_{n}$ is defined as 
  \be
     |M|^2 = {M\ast\widetilde{M}}= {\sum_{A} M_A^2}.
  \ee
  For $n=2(\mod 4)$ and $n=3(\mod 4)$ the \blue{pseudoscalar} is
  $i_n=\vect{e}_1\vect{e}_2\ldots\vect{e}_n$ with 
  (also valid for $Cl_{0,n}, \,n=1,2(\mod 4)$)
  \begin{equation}
     \alert{i_n^2=-1}.
  \end{equation} 
  A \blue{blade} $B$ describes a vector \blue{subspace} 
  \be
    V_B=\{ \vect{x}\in \R^{p,q} | \vect{x}\wedge B =0 \}.
  \ee 
  Its \blue{dual blade} 
  \be
  B^{\ast}= Bi_n^{-1}
  \ee
  describes the \blue{complimentary} vector subspace $V^{\perp}_B$. 
  The pseudoscalar $i_n \in Cl_n$ is \alert{central} for $n=3(\mod 4)$
  \be
  \qquad i_n \, M = M \, i_n , \qquad \forall M \in Cl_{n}.
  \ee
But for even $n$ we get due to non-commutativity \cite{HM:ICCA7}
of the pseudoscalar $i_n \in Cl_n $ 
  \begin{equation} \label{eq:evencom}
     i_n  M = M_{even}  \,i_n - M_{odd} \,i_n\,,  
     \quad
     e^{i_n\lambda} M 
     = M_{even} \;e^{i_n\lambda} + M_{odd} \;e^{-i_n\lambda},
     \qquad \forall M \in Cl_n ,  \quad \lambda \in \R.
  \end{equation}
\blue{A multivector valued function} 
  $f: \R^{p,q} \rightarrow Cl_{p,q}, \,\, p+q=n, $ has $2^n$ blade components
  $(f_A: \R^{p,q} \rightarrow \R)$
  \begin{equation}\label{eq:MVfunc}
    f(\mbox{\textbf{\textit{x}}})  =  \sum_{A} f_{A}(\vect{x}) {\vect{e}}_{A}.
  \end{equation}
We define the \blue{inner product} of 
$\R^n \rightarrow Cl_{n}$ functions  $f, g $ by
\begin{equation}\label{eq:mc2}
  (f,g) 
  = \int_{\R^n}f(\vect{x})
    \widetilde{g(\vect{x})}\;d^n\vect{x}
  = \sum_{A,B}\vect{e}_A \widetilde{\vect{e}_B}
    \int_{\R^n}f_A (\vect{x})
    g_B (\vect{x})\;d^n\vect{x},
\end{equation}
and the $L^2(\mathbb{R}^n;Cl_{n})$-\blue{norm} 
\begin{equation}\label{eq:0mc2}
  \|f\|^2 
   = \left\langle ( f,f ) \right\rangle
   = \int_{\mathbb{R}^n} |f(\vect{x})|^2 d^n\vect{x},
   \qquad
   L^2(\R^n;Cl_{n})
   = \{f: \R^n \rightarrow Cl_{n} \mid \|f\| < \infty \}. 
\end{equation}

For the Clifford geometric algebra Fourier transformation (CFT) \cite{HM:ICCA7}
the \blue{complex unit $i \in \C$} is replaced by  
some geometric (square) root of $-1$, e.g. pseudoscalars $i_n$, $n=2,3(\mod 4)$. 
\blue{Complex functions $f$} are replaced by multivector functions $f \in L^2(\R^{n};Cl_{n})$. 
\begin{definition}[Clifford geometric algebra Fourier transformation (CFT)]  
\label{df:CFT}
The Clifford GA Fourier transform 
  $\mathcal{F} \{f\}$: $\R^n \rightarrow Cl_n, \,n=\alert{2},{3}(\rm mod\,4)$ 
is given by
\begin{equation}\label{eqmk1}
  \mathcal{F}\{f\}(\bomega)
  = \widehat{f}(\bomega)
  = \int_{\R^n} f(\vect{x})
    \,e^{-i_n \bomega \cdot \svect{x}}\, d^n\vect{x},
\end{equation}
for multivector functions $f$: $\R^n \rightarrow Cl_n $. 
\end{definition}
NB: The CFT can also be defined analogously for $Cl_{0,n'}, \,n'=\alert{1},\blue{2}(\mod 4)$.
The CFT \eqref{eqmk1} is inverted by 
\begin{equation}\label{eq11}
  f(\vect{x})
  = \mathcal{F}^{-1}[\mathcal{F}\{f\}(\bomega )]
  = \frac{1}{(2\pi)^n} \int_{\R^n}\mathcal{F}\{f\}(\bomega ) \, 
    e^{i_n\bomega \cdot \svect{x}}\, d^n \bomega .
\end{equation}

The \textit{similitude group} $\mathcal{G}=SIM(n)$ of \alert{dilations, rotations} and \alert{translations} is a subgroup of  the affine group of $\mathbb{R}^n$ 
\begin{equation}\label{eq:mc3}
  \mathcal{G} 
  = \mathbb{R}^+ \times SO(n)\rtimes \mathbb{R}^n 
  = \{(a,r_{\sbtheta}, \vect{b})|
    a \in \mathbb{R}^+,r_{\sbtheta} \in SO(n), 
    \vect{b} \in \mathbb{R}^n \}.
\end{equation}
The left Haar measure on $\mathcal{G}$ is given by
\begin{equation}
  d\lambda 
  = d\lambda(a,\btheta,\vect{b}) 
  = d\mu(a,\btheta) d^n\vect{b},
  \quad
  d\mu 
  = d\mu(a,\btheta) 
  = \frac{dad\btheta} {a^{n+1}},
\end{equation}
where $d\btheta$
is the Haar measure on $SO(n)$.
We define the \blue{inner product} of 
$ f,g : \mathcal{G} \rightarrow Cl_{n}$ by
\begin{equation}
  ( f,g ) =
  \int_{\mathcal{G}} f(a, \btheta, \vect{b})
  \widetilde{g(a, \btheta, \vect{b})}\;
  d\lambda(a,\btheta,\vect{b} ),
\end{equation}
and the $L^2(\mathcal{G};Cl_{n})$-\blue{norm} 
\begin{equation}
  \|f\|^2 
  = \left\langle ( f,f )\right\rangle
  = \int_{\mathcal{G}} |f(a,\btheta, \vect{b})|^2 
    d\lambda, 
  \quad
  L^2(\mathcal{G};Cl_{n})
  = \{f: \mathcal{G} \rightarrow Cl_{n} \mid \|f\| < \infty \}. 
\nonumber
\end{equation}

\section{Clifford GA wavelets}

Previous works on Clifford wavelets include \cite{MM:CMRA,Brackx:etal,LT:QW,EBC:pubs,JL:QaAW,KCF:MonW,SB:Wavelets,MH:CliffWUP}.
We represent the transformation group {$\mathcal{G}=SIM(n)$}
by applying
\alert{translations}, \alert{scaling and rotations} to a 
so-called
\textit{Clifford mother wavelet}
$\psi : \R^n \rightarrow Cl_n$
\begin{equation}
\psi(\vect{x}) \longmapsto 
\mbox{$\psi$}_{a,\btheta,\svect{b}}
(\vect{x})
= \frac{1}{a^{n/2}} \psi 
(r_{\btheta}^{-1}(\frac{\vect{x}-\vect{b}}{a})).
\end{equation}
The family of wavelets 
$\psi_{a,\btheta,\svect{b}}$
are so-called \alert{Clifford daughter wavelets}.
\begin{lemma}[Norm identity]
The factor ${a^{{-n}/{2}}}$ in $\psi_{ a,\btheta,\svect{b}}$ ensures (independent of $a, \btheta, \vect{b}$) that
\begin{equation}\label{eq:0mc7}
\|\psi_{ a,\btheta,\svect{b}}\|_{L^2(\mathbb{R}^n;Cl_{n})} 
= \| \psi \|_{L^2(\mathbb{R}^n;Cl_{n})}.
\end{equation}
\end{lemma}
The CFT \textit{spectral representation} of Clifford daughter wavelets is 
\begin{equation}\label{eq:mc7}
  \mathcal{F}\{\mbox{$\psi$}_{ a,\btheta,\vect{b}}\}
   (\bomega) 
  =
  a^{\frac{n}{2}}\widehat{\psi}
  (a r_{\btheta}^{-1}(\bomega))
  e^{-i_n\svect{b}\cdot \bomega} \,.
\end{equation}
A Clifford mother wavelet $\psi \in L^2(\R^n;Cl_{n})$
is \textit{admissible} if
\begin{equation}
C_{\psi} 
  =  \int_{\mathbb{R}^+}\int_{S0(n)}a^n \{\widehat{\psi}
  (a r_{\sbtheta}^{-1}(\bomega))\}^{\sim}
  \widehat{\psi}(a r_{\sbtheta}^{-1}(\bomega))\;d\mu
  = \int_{\mathbb{R}^n} 
  \frac{\widetilde{\widehat{\psi}}(\bomega) \widehat{\psi}(\bomega)}
  {|\bomega|^n} \; d^n\bomega, 
\end{equation}
is an \textit{invertible multivector constant and finite} at  a.e. 
$\bomega \in \R^n$.
We must therefore have $\widehat{\psi}(\bomega=0)=0$ and therefore \alert{every} Clifford mother wavelet \alert{component} 
\begin{equation}
  \int_{\mathbb{R}^n} \psi_A (\vect{x}) \;d^n\vect{x}  
  = 0.
\end{equation}
By construction $C_{\psi} = \widetilde{C_{\psi}}$. Hence for $n=2,3 (\mod 4)$
\begin{equation*}
  \langle C_{\psi} \rangle_0 > 0, 
  \quad
  C_{\psi} = \sum_{k=0}^{[n/4]}(\langle C_{\psi} \rangle_{4k} +  \langle C_{\psi} \rangle_{4k+1}).
\end{equation*}
The invertibility of $C_{\psi}$ depends on its grade content, e.g. for $n=2,3$, $C_{\psi}$ is
\alert{invertible}, if and only if 
$\langle C_{\psi} \rangle_1^2 \neq \langle C_{\psi} \rangle_0^2 \;$:
\begin{equation}\label{eq:ICl}
C_{\psi}^{-1} = \frac{\langle C_{\psi} \rangle_0 -\langle C_{\psi} \rangle_1 }
{{\langle C_{\psi} \rangle_0^2 -\langle C_{\psi} \rangle_1^2 } }.
\end{equation}

\begin{definition}[Clifford GA wavelet transformation]
\label{df:GAWT}
For an admissible GA mother wavelet $\psi \in L^2(\mathbb{R}^n;Cl_{n})$ 
and a multivector signal function
$f \in L^2(\mathbb{R}^n;Cl_{n})$
\begin{eqnarray}\label{eqq4}
  T_{\psi} : L^2(\mathbb{R}^n;Cl_{n}) &\rightarrow & L^2(\mathcal{G};Cl_{n}),
  \qquad  
  f \mapsto T_{\psi}f(a,\btheta,\vect{b})
    = \int_{\mathbb{R}^n}f(\vect{x})
    \widetilde{\psi_{a,\btheta,\vect{b}}(\vect{x})}
    \;d^n\vect{x}.
\end{eqnarray}
\end{definition}
\noindent
\textbf{NB:} Because of \eqref{eq:evencom}
we need to \textit{restrict} the mother wavelet $\psi$ for $n=2 (\mod 4)$ to even or odd grades: 
Either we have a \textit{spinor wavelet} 
$\psi \in L^2(\mathbb{R}^n;Cl^+_{n})$ with $\varepsilon = 1$,  
or we have an \textit{odd parity vector wavelet}
$\psi \in L^2(\mathbb{R}^n;Cl^-_{n})$ with $\varepsilon = -1$.
NB:  
For $n=3 (\mod 4) $, no grade restrictions exist. We then always have $\varepsilon=1$.
We immediately see from Definition \ref{df:GAWT} that the Clifford GA wavelet transform 
is \textit{left linear} with respect to multivector constants $ \lambda_1,\lambda_2 \in Cl_{n}$.

The {\textit{spectral} (CFT) representation}
of the Clifford wavelet transform is
\begin{equation}\label{eq:cwf}
\hspace*{-1mm}T_{\psi}f(a,\btheta,\vect{b})
  = \frac{1}{(2\pi)^n} 
  \hspace*{-1mm}\int_{\mathbb{R}^n}\hspace*{-0.5mm}
  \widehat{f}(\bomega)\,
  a^{\frac{n}{2}} 
  \{\widehat{\psi}(a r_{\btheta}^{-1}(\bomega))\}^{\sim}
  e^{\varepsilon i_n \svect{b}\cdot \bomega} \, d^n\bomega.
\end{equation}
\textbf{NB}: The CFT for $n=2(\mod 4)$ preserves even and odd grades.

We further have the following set of properties. \textit{Translation covariance}: 
If the argument of $T_{\psi}f(\vect{x})$
is \alert{translated} by a constant $\vect{x}_0 \in \mathbb{R}^n$ then
\begin{equation} 
  [T_{\psi}f(\cdot -\vect{x}_0)](a,\btheta, \vect{b}) 
   = T_{\psi}f(a, \btheta,\vect{b}-\vect{x}_0)\,.
\end{equation}
\textit{Dilation covariance}:
If $0<c\in \R$ then
\begin{equation} 
  [T_{\psi}f(c\,\cdot)](a,\btheta,\vect{b}) 
  = \frac {1}{c^{\frac{n}{2}}} T_{\psi}f(ca, \btheta, c\vect{b})\;.
\end{equation}
\textit{Rotation covariance}:
If  $ r   = r_{\sbtheta} $, 
    $ r_0 = r_{\sbthetaz} $
and $ r'  = r_{\sbthetap} = r_0 r = r_{\sbthetaz} r_{\sbtheta} $
are \alert{rotations}, then
\begin{equation}\label{eq:rotC} 
  [T_{\psi}f(r_{\sbthetaz}\cdot)] (a,\btheta,\vect{b}) 
  = T_{\psi}f(a, \btheta^{\prime}, r_{\sbthetaz}\vect{b})\,.
\end{equation}

Now we see some differences from the classical wavelet transforms. The next property is an 
\textit{inner product relation}:
Let $C^{\,\prime}_{\psi}=(-\varepsilon)^nC_{\psi}$, and 
$f, g \in L^2(\mathbb{R}^n;Cl_{n})$ arbitrary. 
Then we have
\begin{eqnarray}\label{eqM:C1}
  ( T_{\psi}f,T_{\psi}g )_{L^2(\mathcal{G};Cl_{n})} 
  & = &
  ( fC^{\,\prime}_{\psi}, g )_{L^2(\mathbb{R}^n;Cl_{n})} \,.
\end{eqnarray}
The spectral representation \eqref{eq:cwf} and the CFT Plancherel theorem \cite{HM:ICCA7}
are essential for the proof. 
As a corollary we get the following \textit{norm relation}:
\begin{eqnarray}\label{eq1:col}
\| T_\psi f \|_{L^2(\mathcal{G};Cl_{n})}^2 
  & = &
  Sc ( f C^{\,\prime}_{\psi} , f )_{L^2(\mathbb{R}^n;Cl_{n})}
  \,\,\, = \,\,\, C^{\,\prime}_{\psi}\ast ( f , f )_{L^2(\mathbb{R}^n;Cl_{n})} \,.  
\end{eqnarray}
We can further derive the
\begin{theorem}[Inverse Clifford $ Cl_{n}$ wavelet transform]
Any $ f \in L^2(\mathbb{R}^n;Cl_{n})$ can be decomposed with respect to
an admissible Clifford GA wavelet as
\begin{align}\label{eq:iivw1}
  f(\vect{x})  
  = \int_{\mathcal{G}}
    T_{\psi}f (a,\btheta,\vect{b})\,
    \psi_{a,\btheta,\svect{b}}\,
    C_{\psi}^{\,\prime\,-1}\,d\mu d^n\vect{b}
  = \int_{\mathcal{G}}
    ( f, \psi_{a,\btheta,\svect{b}}
    )_{L^2(\mathbb{R}^n;Cl_{n})}
    \psi_{a,\btheta,\svect{b}}
    C_{\psi}^{\,\prime\,-1}\,d\mu d^n\vect{b} ,
\end{align}
the integral converging in the weak sense.
\end{theorem}
Next is the \textit{reproducing kernel}:
We define for an admissible Clifford mother wavelet $\psi \in L^2(\mathbb{R}^n;Cl_{n})$
\begin{equation}
  \BK_{\psi}(a,\btheta,\vect{b};
  a^{\prime},\btheta^{\prime},\vect{b}^{\prime}) 
  = 
  \left( 
    \psi_{a,\btheta,\svect{b}}
    C_{\psi}^{\,\prime\,-1} , 
    \psi_{a^{\prime},\btheta^{\prime},\svect{b}^{\prime}}
  \right)_{L^2(\mathbb{R}^n;Cl_{n})} \,.
\end{equation}
Then 
$\BK_{\psi}(a,\btheta,\vect{b};a^{\prime},\btheta^{\prime},
\vect{b}^{\prime})$
is a reproducing kernel in 
$L^2(\mathcal{G}, d\lambda )$,
i.e,
\begin{equation}
  T_{\psi} f (a^{\prime}, \btheta^{\prime},\vect{b}^{\prime}) 
  =
  \int_{\mathcal{G}}T_{\psi}f (a,\btheta,\vect{b})
  \BK_{\psi}(a,\btheta,\vect{b};a^{\prime},\btheta^{\prime},\vect{b}^{\prime}) 
  d\lambda \,.
\end{equation}
\begin{theorem}[Generalized GA wavelet uncertainty principle]
\label{th4.1}
Let $\psi $ be an admissible Clifford algebra mother wavelet.
Then for every $f \in L^2(\mathbb{R}^n;Cl_{n})$, 
the following inequality holds 
\begin{align}
  \label{eq:genCWUP}
   \| \btheta T_{\psi}f(a,\btheta,\vect{b}) \|^2_{L^2({\mathcal G};Cl_{n})}
   \;\;
   C_{\psi} \ast 
   &(\widetilde{\bomega\hat{f}} , \widetilde{\bomega \hat{f}})_{L^2(\mathbb{R}^n;Cl_{n})}
   \geq  
   \frac{n(2\pi)^n}{4}\,\,\left[C_{\psi}\ast ( f,f )_{L^2(\mathbb{R}^n;Cl_{n})}\right]^2.
\end{align}
\end{theorem} 
\noindent
\textbf{NB:} The integrated variance
$
  \int_{\R^+}\int_{SO(n)}
  \|\bomega \mathcal{F}\{T_{\psi}f(a,\btheta, . ) \} \|^2_{L^2(\mathbb{R}^n;Cl_{n})} d \mu
$
is independent of the wavelet parity $\varepsilon$. Otherwise the {proof} is similar to the one for $n=3$ in \cite{MH:CliffWUP}. For scalar admissibility constant this reduces to
\begin{corollary}[Uncertainty principle for GA wavelet]
  \label{cor:scUP}
Let $\psi $ be a Clifford algebra  wavelet with scalar 
admissibility constant.
Then for every $f \in L^2(\mathbb{R}^n;Cl_{n})$, the following inequality holds 
\begin{equation}\label{eq:UPCWTs}
  \|\vect{b} \,T_{\psi} f(a,\btheta,\vect{b})\|_{L^2(\mathcal{G};Cl_{n})}^2\;
  \|\bomega \hat{f}\|_{L^2(\mathbb{R}^n;Cl_{n})}^2
  \geq 
  n C_{\psi} \frac{(2\pi)^n}{4} \|f\|^4_{L^2(\mathbb{R}^n;Cl_{n})}.
\end{equation}
\end{corollary}

Finally \textit{GA Gabor Wavelets} are defined as 
(variances $\sigma_k, 1 \leq k \leq n$, for
$n=2(\mod 4): A \in Cl^+_n$ or $A\in Cl^-_n$)
\begin{equation}
  \psi^{c}(\vect{x}) 
  = 
  \frac{A}{(2\pi)^{\frac{n}{2}} \prod_{k=1}^n \sigma_k }\,
  e^{-\frac{1}{2}\sum_k\frac{x_k^2}{\sigma_k^2} }
\left(
  e^{i_n \bomega_0\cdot \svect{x}}
  - \underbrace{e^{-\frac{1}{2}\sum_{k=1}^n \sigma_k^2 \omega_{0,k}^2}}_{\mbox{constant}}
\right),
\quad \vect{x},\bomega_0 \in \R^n, \quad \mbox{constant }A \in Cl_n\,.
\end{equation}

\begin{theacknowledgments}
  Soli deo gloria. I do thank my dear family, B. Mawardi, G. Sommer, W. Spr\"{o}ssig and K. G\"{u}rlebeck.
\end{theacknowledgments}

\bibliographystyle{aipproc}

\end{document}